\documentclass[english,12pt]{smfart}
\usepackage{amssymb}
\usepackage{amsfonts}
\usepackage{amscd}
\usepackage[all]{xypic}

\usepackage[T1]{fontenc}

\CompileMatrices

\swapnumbers

\def\ord{{\rm ord}}
\def\ac{{\overline{\rm ac}}}

\def\11{{\mathbf 1}}
\def\AA{{\mathbf A}}

\def\NN{{\mathbf N}}

\def\QQ{{\mathbf Q}}
\def\RR{{\mathbf R}}

\def\ZZ{{\mathbf Z}}

\def\cL{{\mathcal L}}
\def\cM{{\mathcal M}}

\def\cO{{\mathcal O}}

\mathchardef\alphag="7C0B \mathchardef\betag="7C0C
\mathchardef\gammag="7C0D \mathchardef\deltag="7C0E
\mathchardef\varepsilong="7C22 \mathchardef\varphig="7C27
\mathchardef\psig="7C20 \mathchardef\zetag="7C10
\mathchardef\epsilong="7C0F \mathchardef\rhog="7C1A
\mathchardef\taug="7C1C \mathchardef\upsilong="7C1D
\mathchardef\iotag="7C13 \mathchardef\thetag="7C12
\mathchardef\pig="7C19 \mathchardef\sigmag="7C1B
\mathchardef\etag="7C11 \mathchardef\omegag="7C21
\mathchardef\kappag="7C14 \mathchardef\lambdag="7C15
\mathchardef\mug="7C16 \mathchardef\xig="7C18
\mathchardef\chig="7C1F \mathchardef\nug="7C17
\mathchardef\varthetag="7C23 \mathchardef\varpig="7C24
\mathchardef\varrhog="7C25 \mathchardef\varsigmag="7C26
\mathchardef\Omegag="7C0A \mathchardef\Thetag="7C02
\mathchardef\Sigmag="7C06 \mathchardef\Deltag="7C01
\mathchardef\Phig="7C08 \mathchardef\Gammag="7C00
\mathchardef\Psig="7C09 \mathchardef\Lambdag="7C03
\mathchardef\Xig="7C04 \mathchardef\Pig="7C05
\mathchardef\Upsilong="7C07

\newtheorem{theorem}[subsection]{Theorem}
\newtheorem{lem}[subsection]{Lemma}
\newtheorem{cor}[subsection]{Corollary}
\newtheorem{prop}[subsection]{Proposition}

\theoremstyle{definition}
\newtheorem{definition}[subsection]{Definition}
\newtheorem{example}[subsection]{Example}

\newtheorem{def-prop}[subsection]{Proposition-Definition}
\newtheorem{def-theorem}[subsection]{Theorem-Definition}
\newtheorem{def-lem}[subsection]{Lemma-Definition}

\theoremstyle{remark}
\newtheorem{remark}[subsection]{Remark}

\theoremstyle{plain}

\numberwithin{equation}{subsection}

\def\boxit#1#2{\setbox1=\hbox{\kern#1{#2}\kern#1}%
\dimen1=\ht1 \advance\dimen1 by #1 \dimen2=\dp1 \advance\dimen2 by
#1
\setbox1=\hbox{\vrule height\dimen1 depth\dimen2\box1\vrule}%
\setbox1=\vbox{\hrule\box1\hrule}%
\advance\dimen1 by .4pt \ht1=\dimen1 \advance\dimen2 by .4pt
\dp1=\dimen2 \box1\relax}

\renewcommand{\theequation}{\thesubsection.\arabic{equation}}

\def\AA{{\mathbf A}}

\def\NN{{\mathbf N}}

\def\QQ{{\mathbf Q}}
\def\RR{{\mathbf R}}

\def\ZZ{{\mathbf Z}}

\def\cL{{\mathcal L}}
\def\cM{{\mathcal M}}

\def\cO{{\mathcal O}}

\mathchardef\alphag="7C0B \mathchardef\betag="7C0C
\mathchardef\gammag="7C0D \mathchardef\deltag="7C0E
\mathchardef\varepsilong="7C22 \mathchardef\varphig="7C27
\mathchardef\psig="7C20 \mathchardef\zetag="7C10
\mathchardef\epsilong="7C0F \mathchardef\rhog="7C1A
\mathchardef\taug="7C1C \mathchardef\upsilong="7C1D
\mathchardef\iotag="7C13 \mathchardef\thetag="7C12
\mathchardef\pig="7C19 \mathchardef\sigmag="7C1B
\mathchardef\etag="7C11 \mathchardef\omegag="7C21
\mathchardef\kappag="7C14 \mathchardef\lambdag="7C15
\mathchardef\mug="7C16 \mathchardef\xig="7C18
\mathchardef\chig="7C1F \mathchardef\nug="7C17
\mathchardef\varthetag="7C23 \mathchardef\varpig="7C24
\mathchardef\varrhog="7C25 \mathchardef\varsigmag="7C26
\mathchardef\Omegag="7C0A \mathchardef\Thetag="7C02
\mathchardef\Sigmag="7C06 \mathchardef\Deltag="7C01
\mathchardef\Phig="7C08 \mathchardef\Gammag="7C00
\mathchardef\Psig="7C09 \mathchardef\Lambdag="7C03
\mathchardef\Xig="7C04 \mathchardef\Pig="7C05
\mathchardef\Upsilong="7C07

\DeclareMathOperator*{\sq}{{\square}}
\def\ord{{\rm ord}}

\begin{document}

\title[Lipschitz continuity over the  $p$-adics]{Lipschitz continuity properties for $p$-adic semi-algebraic and subanalytic functions}

\author{Raf Cluckers}
\address{Katholieke Universiteit Leuven, Department of Mathematics,
Celestijnenlaan 200B, B-3001 Leu\-ven, Bel\-gium\\ The author is a
postdoctoral fellow of the Fund for Scientific Research - Flanders
(Belgium) (F.W.O.)} \email{raf.cluckers@wis.kuleuven.be}
\urladdr{http://www.wis.kuleuven.be/algebra/Raf/}

\author{Georges Comte}

\address{Laboratoire J.-A. Dieudonn\'e,
Universit\'e de Nice - Sophia Antipolis, Parc Valrose, 06108 Nice
Cedex 02, France (UMR 6621 du CNRS)} \email{comte@math.unice.fr}
\urladdr{http://www-math.unice.fr/membres/comte.html}

\author{Fran\c cois Loeser}

\address{{\'E}cole Normale Sup{\'e}rieure,
D{\'e}partement de math{\'e}matiques et applications, 45 rue d'Ulm,
75230 Paris Cedex 05, France (UMR 8553 du CNRS)}
\email{Francois.Loeser@ens.fr}
\urladdr{http://www.dma.ens.fr/$\sim$loeser/}

\begin{abstract}
We prove that a (globally) subanalytic function $f:X\subset
\QQ_p^n\to \QQ_p$ which is locally Lipschitz continuous with
some constant $C$ is piecewise (globally on
each piece) Lipschitz continuous with possibly some
other constant, where the pieces can be taken subanalytic. We also
prove the analogous result for a subanalytic family of functions
$f_y:X_y\subset \QQ_p^n\to \QQ_p$ depending on $p$-adic parameters. The statements also hold in a
semi-algebraic set-up and also in a finite field extension of
$\QQ_p$. These results are
$p$-adic analogues of
results of K.~Kurdyka over the real numbers.
To encompass the total disconnectedness of $p$-adic fields,
we need to introduce new  methods adapted to the 
$p$-adic situation.
\end{abstract}

\maketitle

\renewcommand{\partname}{}

\section*{Introduction}

In the real setting, a $C^1$-function on an interval in $\RR$ which
has bounded derivative is automatically Lipschitz continuous.
Indeed, if $f:(a,b)\to \RR$ satisfies $|f'(x)|\leq C$ for some $C$ and
all $x$, then, for any $c<d$ in $(a,b)$ one has $|f(c)-f(d)|=|
\int_c^d f'(x)dx |\leq C|c-d|$.
Such a result cannot hold
for general $C^1$-functions over the $p$-adics,
because of
total disconnectedness.
Indeed, there are easy examples
of locally constant (hence $C^1$) functions
$g:X\subset \QQ_p\to\QQ_p$ with $X$ open, for which there exists no partition of
$X$  into finitely many pieces such that $g$ is Lipschitz continuous on each piece
(see Example \ref{exloc}).
Also, there are examples
of $C^1$ functions $g: X \subset \QQ_p\to\QQ_p$, with $g'$ identically zero, that are not locally Lipschitz continuous for
any constant $C$, at an infinite number of points (Example \ref{exloc2}).
Such examples show that over $p$-adic fields, the relation between bounds on derivatives and local Lipschitz properties may be quite chaotic,
so
in order to be able
to obtain
 significant
results related to piecewise Lipschitz continuity, it seems reasonable to limit the class of  $p$-adic functions we consider
   to a class
of tame  piecewise $C^1$-functions where in particular functions of this class can be  described by finite amounts
of data.
  In this paper we consider
  two such tame classes of $p$-adic sets and functions:
semi-algebraic sets and
functions on the one hand, and
(globally) subanalytic $p$-adic sets and
functions on the other hand. Semi-algebraic functions are a natural generalization of algebraic functions on algebraic subsets of $\QQ_p^n$, and subanalytic functions are a further enrichment of semi-algebraic functions with certain analytic functions. In both frameworks, the functions are more than piecewise $C^1$: they are even piecewise (locally) $\QQ_p$-analytic where the pieces are $\QQ_p$-analytic manifolds which are moreover subanalytic (resp.~semi-algebraic). In these two frameworks we obtain several results about piecewise
Lipschitz continuity for multi-variable functions $g:X\subset
\QQ_p^m\to \QQ_p$, assuming only local conditions on $g$, like for example boundedness of the partial derivatives.
Of course, a (locally) $\QQ_p$-analytic function $f:X\to \QQ_p$ on an open
set
$X\subset \QQ_p^m$ satisfies the following local property, cf.~Lemma \ref{localeL}:
\begin{quote}
if $|\partial f(x)/\partial x_i|\leq 1$ for all $i=1,\ldots,m$ and all $x$, then $f$ is locally Lipschitz continuous with constant $1$.
\end{quote}
From this one deduces, cf.~Proposition \ref{localL}:
\begin{quote}
 Let $g:X\subset \QQ_p^m\to\QQ_p$ be subanalytic (resp.~semi-algebraic) and $C^1$ on an open $X$ such that $|\partial g(x)/\partial x_i|\leq 1$ for all $i=1,\ldots,m$. Take any finite partition of $X$ into subanalytic (resp.~semi-algebraic) $\QQ_p$-analytic manifolds $X_i$ on which $g$ is $\QQ_p$-analytic. Then the restriction of $g$ to $X_i$ is locally Lipschitz continuous with constant $1$ for each $i$.
\end{quote}

Note that the manifolds $X_i$ are not necessarily open in $\QQ_p^m$.
This indicates  that, for functions whose domain is not necessarily open, it may show more convenient to work with
the condition of local Lipschitz continuity instead of conditions on partial derivatives.

Over the reals Kurdyka \cite{Kurdyka} obtains the following result: if
a (globally) subanalytic function $f:X\subset \RR^n\to \RR$ is
locally Lipschitz continuous for some fixed constant $C$, then it is piecewise (globally on
the piece) Lipschitz continuous with possibly some other constant.
 We prove the following $p$-adic analogue: if a subanalytic
function $g:X\subset \QQ_p^n\to \QQ_p$ is locally Lipschitz
continuous with constant $C$, then it is piecewise (globally on the piece) Lipschitz
continuous for possibly some other constant.
 In both the real and the $p$-adic setting, the pieces can be taken
to be subanalytic, resp.~semi-algebraic if $f$ and $g$ are
semi-algebraic. In fact, we will prove this result for a fixed finite field extension  $K$ of $\QQ_p$, and for subanalytic families of functions instead of for individual subanalytic functions.

\subsection*{}
Let us start by  explaining the real (globally) subanalytic and
semi-algebraic situation, according to  Kurdyka \cite{Kurdyka},
giving a rough
sketch of
 the main arguments for functions in
up to two variables. The one variable case becomes trivial by the
mentioned relation between distance and the integral of the
derivative once one notes that any globally subanalytic subset of
$\RR$ is a finite union of points and open intervals (compare with
$o$-minimality). For $f:X\subset \RR^2\to\RR$ a globally subanalytic
function, Kurdyka proves that $X$ can be partitioned into finitely
many pieces which are $o$-minimal cells
(sometimes
called cylinders),
 for example of the form
$$
A_1=\{(x,y)\in(a,b)\times \RR \mid \alpha(x) <y<  \beta(x) \},
$$
or of the form
$$
A_2=\{(x,y)\in(a,b)\times \RR \mid \alpha(x) =y\},
$$
 where moreover the ``boundary functions'' 
 $\alpha$ and $\beta$ 
 are (globally) subanalytic or semi-algebraic and
have bounded derivatives. For such a decomposition in
so-called $L$-regular cells to exist, it is important that for each
piece $A$ separately, affine coordinates on $\RR^2$ are chosen
so that $A$ has indeed such a nice description (see also \cite{Pawlucki}).
 Using such a decomposition into
$L$-regular cells and affine coordinates adapted to each cell individually, the Lipschitz continuity result follows using a
path integral of the derivative of $f$ along a well chosen path
inside the cell, after noticing that any two points in an $L$-regular cell
can be connected by a path whose length is not much bigger than the
distance between the two points it connects.

\subsection*{}
Over $\QQ_p$, it seems not possible to follow a strategy similar to Kurdyka's, since
there is no clear notion of paths connecting two points,
let alone of the length of a path,
which are basic ingredients for Kurdyka's approach over the reals.
More generally speaking, as far
as we know there is no clear connection between integrals and
distances between points in $\QQ_p$. Hence, a new approach had to be
devised.
Let us sketch this new approach in the one and the two
variable case. Let $g:X\subset\QQ_p\to\QQ_p$ be a subanalytic
function (globally subanalytic, as always in this paper) which is
locally Lipschitz continuous with constant $C$. We know that we can partition $X$ into
finitely many $p$-adic cells, all of which are very roughly of a
form similar to
 \begin{equation*}
A=\{t\in K\mid |\alpha| \leq |t-c| \leq |\beta|,\
  t-c\in \lambda Q_{m,n}\},
\end{equation*}
with constants $n>0$, $m>0$,
$\lambda\in K$,
and where $Q_{m,n}$ is the set of all
$p$-adic numbers of the form $p^{na}(1+p^{m}x)$ for some $x\in \ZZ_p$ and some
$a\in\ZZ$ (see \ref{def::cell} and \ref{thm:CellDecomp} for precise
statements). We call $c$ the center of $A$ (note that $c$
 may
 lie
outside $A$, which happens precisely when $\lambda\not=0$). We define the
balls of the cell $A$ as the collection of maximal balls (with
respect to inclusion) contained in $A$ (cf.~\ref{ballcell}). It then
follows from a certain Jacobian property \ref{jacprop} that we can
select
 the cells $A$ in such a way that each ball of $A$ is mapped to
either a point or a ball under $g$.
We refine this Jacobian property so that we
can ensure that the images of the balls of $A$ form up to a single
cell which has moreover as collection of maximal balls precisely the collection of the images of the maximal balls in $A$ (cf.~Proposition
\ref{celcel}).
In particular, $g(A)$ is a cell with a certain center
$d$. Then, roughly, distances between points in $A$ are compared to distances to $c$ and similarly in the range of $g$, using the center
$d$. A calculation based on this comparison of distances and the Jacobian property  then finishes
the one variable case. Actually, these arguments also show that the
analogue
statement
for subanalytic \emph{families} of
functions $g_y:X_y\subset \QQ_p\to\QQ_p$ (instead of for individual subanalytic functions) holds. Such families are
used in the several variable case for the induction argument.

In the two variables case,
we obtain a result
vaguely reminiscent of
 real partitions into $L$-regular cells. Let $g:X\subset
\QQ_p^2\to \QQ_p$ be locally Lipschitz continuous with some constant $C$. Roughly, we partition the
family $X_{x_1}:=\{x_2\mid x\in X\}$ into finitely many families of cells
$A_{x_1}\subset K$
with center $c$ and  boundaries $\alpha$ and
$\beta$ now depending on $x_1$. We show that, after possibly
switching
 the role of $x_1$ and $x_2$, we can ensure that the center $c$ is
Lipschitz continuous in $x_1$ (see Proposition \ref{Lipcenter}). By a piecewise
bi-Lipschitz transformation, we may then assume that the center is
identically zero for each of the cells. This is already an important reduction, but
the obstacles
due to the lack of a good notion of paths and integrals to control distances remain.
Instead of
working with paths as is done in the real case, we work with a
finite sequence of points with given
starting
 point and endpoint,
and one could understand such a finite sequence of jumps from one point to the
next as a $p$-adic analogue of a real ``path''. For such a sequence
of jumps
to be of use, the following is required:
after each jump, one
should still stay in the same cell so that one can still evaluate
the function $g$, the total (cumulative) distance of the jumps should be
comparable to the distance between the
starting
 point and the
endpoint, and the function should not vary too much at each jump so
that one can control $|g(a)-g(b)| $ for any jump from $a$ to $b$ in the sequence.
 This is done in the two variable case as follows.
 Let $a$ and $a'$ be given in $A$. Either $\alpha(x_1)$ has bounded
derivative, and then we can use induction for the one-variable
function $x_1\mapsto g(x_1,\alpha(x_1))$ and roughly jump from $a$
to $(a_1,\alpha(a_1)))$, then to $(a_1',\alpha(a_1')))$ and finally
to $a'$.
  In the at first sight more difficult case
  where $\alpha(x_1)$ has large
derivative, we invert the role of $x_1$ and $x_2$ in the
parametrization of the function $f(x_1,\alpha(x_1))$, namely, we
essentially
work with the one-variable function $x_2\mapsto f(\alpha^{-1}(x_2),x_2)$ and use
induction for this function and then roughly make similar jumps as
before: from $a$ to $(a_1,\alpha(a_1))) = (\alpha^{-1}(b),b) $, then
to $(\alpha^{-1}(b'),b')=(a'_1,\alpha(a'_1)))$ and finally to $a'$,
for some $b,b'\in \QQ_p$. This ``path'' allow us to bound
$|g(a)-g(a')|$ in terms of $|a-a'| $ as needed for Lipschitz
continuity, uniformly in $a$ and $a'$ in the cell. Of course,
some fine tuning    is required in order to guarantee
 injectivity  before one starts inverting functions like $\alpha$, which is provided e.g. by Corollary \ref{inj-cons}.
 As already indicated, the bulk of the paper is concerned with a fixed finite field extension  $K$ of $\QQ_p$.
We conclude the paper by extending our main theorem to elementary extensions of $K$, see Proposition \ref{ext}.

\subsection*{} This paper arose from our work
  \cite{CCL} which
provides $p$-adic analogues of  real results in \cite{KR},
\cite{Comte}, \cite{Comte2},
(see
also
\cite{ComteMerle} for a multidimensional version),
and of complex results in \cite{Thie}.
The main result of the present paper,
Theorem \ref{lip},  is used in \cite{CCL} to prove the existence of distinguished tangent
cones of definable sets and to establish the $p$-adic
counterpart of Thie's formula of \cite{Thie}.

\subsection*{}

{During the preparation of this paper, the authors have been partially supported by grant
ANR-06-BLAN-0183.}

\section{Basic terminology and results}

\subsection{}
Let $K$ be a fixed finite field extension of $\QQ_p$, the field of
$p$-adic numbers. Write $\cO_K$ for the valuation ring of $K$ and
$\cM_K$ for the maximal ideal of $\cO_K$. We denote by
$\ord:K^\times\to \ZZ$ the valuation and we set $\vert x \vert := q^{-
\ord (x)}$ and $\vert 0 \vert = 0$, with $q$ the cardinality of the
residue field of $K$. For a tuple $x\in K^n$, we write $|x|$ for
$\max_{i=1}^n |x_i|$.

A \emph{ball} in $K$ is a subset of the form $a+b\cO_K$ with $a\in
K$ and $b$ in $K^\times$. Note that in this terminology, a ball is always a nonempty, open, proper subset of $K$.

Let $\pi_K$ be a uniformizer of $\cO_K$. For each integer $n>0$, let
$\ac_n:K\to \cO_K/(\pi_K^n)$ be the map sending $0$ to $0$ and
nonzero $x$ to $x \pi_K^{-\ord(x)}\bmod (\pi_K^n)$.

The language $\cL$ is, consequently in the whole paper, either the
subanalytic (as in e.g.~\cite{DHM}) or the semi-algebraic language
on $K$ (Macintyre's language), with coefficients (also called parameters) from
$K$. Hence, $\cL$-definable means either subanalytic or
semi-algebraic with parameters from $K$ consequently throughout
the paper, which corresponds to the set-up of \cite{CCL}.

\begin{definition}\label{def:lip}
Given two metric spaces $(X, d_X)$ and $(Y, d_Y)$, where $d_X$
denotes the metric on the set $X$ and $d_Y$ the metric on
$Y$, a function $f: X \to Y$ is called Lipschitz continuous if there
exists a real constant $C \geq  0$ such that, for all $x_1$ and
$x_2$ in $X$,
$$
    d_Y(f(x_1), f(x_2)) \leq C d_X(x_1, x_2).
$$
In the above case, we also call $f$
Lipschitz continuous with constant $C$, or just $C$-Lipschitz continuous. If there is a constant $C$
such that locally around each $x\in X$ the function $f$ is
$C$-Lipschitz continuous, then $f$ is called locally Lipschitz
continuous with constant $C$, or just locally $C$-Lipschitz
continuous.
\end{definition}
 In this paper, the metrics come from the $p$-adic norm on the spaces $K^n$.
 For general $K$-analytic functions (which are not necessarily $\cL$-definable) on an open domain in $K^m$, there is a general link between bounded partial derivatives and  local Lipschitz continuity, as follows.
\begin{lem}[Lemma 1.4.6 of \cite{CCL}]\label{localeL}
Let $X\subset K^m$ be open and let $f:X\to K$ be $K$-analytic, meaning that locally, $f$ is given by converging power series over $K$. Suppose that
$$|\partial f(x)/\partial x_i|\leq 1$$
for all $i=1,\ldots,m$ and all $x$ in $X$. Then $f$ is locally Lipschitz continuous with constant $1$.
\end{lem}

If $f$ is merely $C^1$ one has to be more careful, even for one variable functions (cf. Example \ref{exloc2}). The next proposition is about $\cL$-definable $C^1$ functions on an open $X$.

\begin{prop}\label{localL}
Let $X\subset K^m$ be open and $\cL$-definable and let $f:X \to K$ be $\cL$-definable. Suppose that $f$ is $C^1$, and that
$$
|\partial f(x) / \partial x_j| \leq  1
$$
for all $x$ in $X$ and for all $j=1,\ldots,m$. Let $\{X_i\}_i$ be any finite partition of $X$ into $\cL$-definable parts which are $K$-analytic manifolds on which $f$ is $K$-analytic (such partition always exists). Then the restriction of $f$ to $X_i$ is locally Lipschitz continuous with constant $1$ for each $i$.

Moreover, the same results hold for $\cL$-definable families of open  $X_y\subset K^m$ and $C^1$ functions $f_y:X_y\to K$, where $y$ runs over an $\cL$-definable set $Y$.
\end{prop}
Note that in Proposition \ref{localL}, the $X_i$ need not be open in $K^m$, hence they are different from the situation of Lemma \ref{localeL}.
\begin{proof}
That such a finite partition exists follows from the Cell Decomposition Theorem \ref{thm:CellDecomp} below, but was already obtained in \cite{DvdD} without cell decomposition. Take a point $x_0$ on some $X_i$. Suppose that the manifold $X_i$ is of dimension $d$. By the implicit function theorem and by the non archimedean property, there exists an open neighborhood $U$ of $x_0$ in $X_i$ and a $K$-bi-analytic isometry $i:U\subset K^m\to B^d$ for some ball $B\subset K$. We can finish by Lemma \ref{localeL} applied to $f\circ i^{-1}$. Exactly the same proof works for families $f_y:X_y\to K$. Indeed, by cell decomposition and up to a finite $\cL$-definable partition of the family $X_y$ we may suppose that $X_y$ is a $K$-analytic manifold on which $f_y$ is $K$-analytic for each $y$.
\end{proof}

\begin{example}\label{exloc2}
Proposition \ref{localL} has no analogue for general $C^1$ functions, even in just one variable, say, from $K$ to $K$, using finite partitions. For example, write $K$ as a countable disjoint union of translates of the ball $\cM_K$,
$$
K= \bigcup_{i\in\NN} a_i + \cM_K
$$
for some choice of the $a_i\in K$.
Let $f:K\to K$ send $a_i+x$ with $x\in \cM_K$ to $g(x)$ where $g:\cM_K\to K$ is defined as follows.
Write $\cM_K\setminus \{0\}$ as a countable union of disjoint balls of the form $b+b^3\cO_K$, that is,
$$
\cM_K \setminus \{0\} = \bigcup _{i\in\NN} b_i+b_i^3\cO_K
$$
for some choice of $b_i\in \cM_K\setminus \{0\}$. For each integer $n>0$, fix one of the $b_i$ with $\ord (b_i)=n$ and call these fixed $b_i$ special.
For $x\in \cO_K$, define $g(b_i+b_i^3x)$ as $0$ if  $b_i$ is non special and as $b_i^2$ if  $b_i$ is special and put $g(0)=0$.
Then $g$ and $f$ are $C^1$ and $f'$ and $g'$ are both identically zero.
Hence, one can take $C$ arbitrarily small. However, $g$ is not locally $C'$-Lipschitz continuous around $0$ for any constant $C'>0$. Indeed, let $B$ be a small enough ball around $0$, and take a special $b_i$ inside $B$ close enough to zero. Then, there exists a non special $b_j$ inside $B$ such that
$$
|b_i-b_j| = q \cdot |b_i^3|,
$$
with $q$ the cardinality of the
residue field of $K$.
On the other hand,
$$
| g(b_i)-g(b_j)   | = | b_i^2 |.
$$
Hence, for $f$, there exists no finite partition of $K$ which makes $f$ locally $C'$-Lipschitz continuous on the pieces, for any choice of $C'$.
\end{example}

\section{The main results}

\begin{theorem}[Main theorem]\label{lip}
Let $\varepsilon>0$ be given.  Let $f:X\subset K^m\to K$ be an
$\cL$-definable function which is locally $\varepsilon$-Lipschitz
continuous. Then there exist $C>0$ and a finite definable partition
of $X$ into parts $A_i$ such that the restriction of $f$ to $A_i$ is
(globally) $C$-Lipschitz continuous for each $i$.
\end{theorem}

\begin{example}\label{exloc}
Theorem \ref{lip} has no analogue for general $C^1$ functions, even in just one variable, say, from an open $X\subset K$ to $K$. For example, if $X$ is the open $\cO_K\setminus \{0\}$, and $f$ sends $x\in X$ to $|x|$, where the rational number $|x|$ is seen as an element of $K$, then $f$ is clearly locally constant, but, for $x_1, x_2 \in X$ with $|x_2|<|x_1|$ one has $|f(x_1)-f(x_2)| \geq |x_2|^{-1}$ which grows to infinity while $|x_1-x_2|=|x_1|$ goes to zero whenever $x_1$ approaches zero in $X$.
\end{example}

In order to formulate a variant of Theorem \ref{lip} for families of functions, the following notation will be convenient.
For $g:D\subset A\times B\to C$ a function, and for $b\in B$,  write $g(\cdot,b)$ for the function which sends $a$ with $(a,b)\in D$ to $g(a,b)$.
The domain of $g(\cdot,b)$ is thus the set $\{a\in A\mid (a,b)\in D\}$ which we will denote by $D_b$.

\begin{theorem}[Main theorem: parameterized version]\label{rellip}
Let $\varepsilon>0$ be given. Let $Y$ be an $\cL$-definable set. Let
$f:X\subset K^m\times Y\to K$ be an $\cL$-definable function such
that for each $y\in Y$ the function $f(\cdot,y):x\mapsto f(x,y)$ is
locally $\varepsilon$-Lipschitz continuous on $X_y$. Then there exist $C>0$
and a finite definable partition of $X$ into parts $A_i$ such that
for each $y\in Y$ and $i$ the restriction of $f(\cdot,y)$ to
$A_{iy}$ is (globally) $C$-Lipschitz
continuous.
\end{theorem}

The following proposition compares to the  notion of $L$-regular
cells on the real number field, see \cite{Kurdyka}, which goes
back to A.~Parusi\'nski \cite{Paru}, see also the more recent \cite{Pawlucki}. The definition of $p$-adic
cells and their centers will be given in section \ref{s:cells}.

\begin{prop}[Cells with Lipschitz continuous centers]\label{Lipcenter}
Let $Y$ and $X\subset K^m\times Y$ be $\cL$-definable. Then there
exist $C>0$, a finite partition of $X$ into $\cL$-definable parts
$A$ and for each part $A$ a coordinate projection
$$
\pi:K^m\times Y\to K^{m-1}\times Y
$$
such that, over $K^{m-1}\times Y$ along this projection $\pi$, the
set $A$ is a $p$-adic cell with center $c:\pi(A) \to K$ and such that moreover the
function
 $$
c(\cdot,y): (x_1,\ldots,x_{m-1})\mapsto
c(x_1,\ldots,x_{m-1},y)
 $$
is $C$-Lipschitz continuous on $\pi(A)_y$ for each $y\in Y$.
\end{prop}

A last, more technical new result in this paper is Proposition
\ref{celcel}, which, in the one-variable case, says that for an
injective definable function $f$ one can partition the domain and
the range compatibly into cells in some strong sense related to the maximal balls contained in the cells.

\section{Some results related to cell decomposition over $K$}\label{s:cells}

For integers $m>0$ and $n>0$, let $Q_{m,n}$ be the ($\cL$-definable) set
$$
Q_{m,n}:=\{x\in K^\times \mid \ord (x) \in n\ZZ,\ \ac_m(x)=1\}.
$$
For $\lambda\in K$ let $\lambda \cdot Q_{m,n}$ denote $\{\lambda
x\mid x\in Q_{m,n}\}$. The sets $Q_{m,n}$ are a variant of Macintyre's predicates $P_\ell$ of $\ell$th powers; the corresponding notions of cells are slightly different but equally powerful and similar in usage. Indeed, any coset of $P_\ell$ is a finite disjoint union of cosets of some $Q_{m,n}$ and vice versa.

 \begin{definition}[$p$-adic cells]\label{def::cell}
Let $Y$ be an $\cL$-definable set.
A $p$-adic cell $A\subset K\times Y$ over $Y$ is a (nonempty) set of
the form
 \begin{equation}
A=\{(t,y)\in K\times Y\mid y\in Y',\ |\alpha(y)| {\sq}_1 |t-c(y)|{\sq}_2
|\beta(y)|,\
  t-c(y)\in \lambda Q_{m,n}\},
\end{equation}
with $Y'$ a $K$-analytic $\cL$-definable manifold, constants $n>0$, $m>0$, $\lambda$
in $K$, $\alpha,\beta:Y'\to K^\times$ and $c:Y'\to K$ all $K$-analytic
$\cL$-definable functions, and $\square_i$ either $<$ or no
condition, such that $A$ projects surjectively onto $Y'\subset Y$.
 We call $c$ the center of the cell $A$, $\lambda Q_{m,n}$ the coset
of $A$, $\alpha$ and $\beta$ the boundaries of $A$, and $Y'$ the
base of $A$. If $\lambda=0$ we call $A$ a $0$-cell
and if $\lambda\not=0$ we call $A$ a $1$-cell.
\end{definition}

Note that a $p$-adic cell over $Y$ is an $\cL$-definable set which is moreover a $K$-analytic manifold.

\begin{def-prop}[Balls of cells]\label{ballcell}
Let $Y$ be $\cL$-definable. Let $A\subset K\times Y$ be a $p$-adic
$1$-cell over $Y$ with coset $\lambda Q_{m,n}$ and center $c$. Then, for each
$(t,y)\in A$ with $y\in Y$, there exists a unique maximal ball
$B_{t,y}$ containing $t$ and satisfying $B_{t,y}\times \{y\}\subset
A$, where the maximality is for the inclusion. We call the
collection of balls $\{B_{t,y}\mid (t,y)\in A\}$ the balls of the cell
$A$; for fixed $y_0\in Y$ we call the collection of balls
$\{B_{t,y_0}\}_{\{t\mid (t,y_0)\in A\}}$ the balls of the cell $A$
above $y_0$. Moreover, for each $(t,y)\in A$ one has
$$
B_{t,y}=\{w\in K\mid \ord (w-c(y)) = a ,\ \ac_m(w-c(y)) =
\ac_m(\lambda)\}
$$
for a unique $a\in \ZZ$ depending on $t$ and $y$.
 If $A\subset K\times Y$ is a $p$-adic $0$-cell then we define the collection of balls of $A$ to be the empty collection, that is, there are no balls of $A$.
\end{def-prop}
\begin{proof}
The uniqueness of $B_{t,y}$ follows from the non archimedean
property. We prove existence of a maximal ball $B_{t,y}$ containing
$t$ and satisfying  $B_{t,y}\times \{y\}\subset A$. Choose $(t,y)\in
A$. Since the collection of balls is preserved under translation by
a constant, we may suppose that $c(y)=0$. Then $(0,y)\not\in A$ since $\lambda\not=0$ and thus $t\not =
0$. Clearly, for $B_1$ being the ball $B_1=t\cO_K$
one has that $B_1\times \{y\}$ is not a subset of $A$ since  $B_1$ contains $0$. One the other hand, let $B_2$ be the ball $t + \pi_K^{m}
t\cO_K$, then clearly $B_2\times \{y\}\subset A$.  Since
the value group is discrete and since $t\in B_2\subset B_1$, the existence follows. In fact, $B_{t,y}=B_2$ since for any strictly bigger ball $B_3$ containing $B_2$ there exists $t'\in B_3$ with $\ac_m(t')\not=\ac_m(t)=\ac_m(\lambda)$. Hence, the description for  $B_{t,y}$ in the proposition follows.
\end{proof}

In the $p$-adic semialgebraic case, Cell Decomposition Theorems are
due to Cohen \cite{cohen} and Denef \cite{D84}, \cite{Dcell} and
they were extended in \cite{Ccell} to the subanalytic setting where
one can find the following version:

\begin{theorem}[$p$-adic Cell Decomposition]\label{thm:CellDecomp}
Let $X\subset K^{m+1}$ and $f_j:X\to K$ be $\cL$-definable for
$j=1,\ldots,r$. Then there exists a finite partition of $X$ into
$p$-adic cells $A_i$ (over $K^m$) with center $c_i$ and coset
$\lambda_i Q_{m_i,n_i}$ such that
 \begin{equation*}
 |f_j(x,t)|=
 |h_{ij}(x)|\cdot|(t-c_i(x))^{a_{ij}}\lambda_i^{-a_{ij}}|^\frac{1}{n_i},\quad
 \mbox{ for each }(x,t)\in A_i,
 \end{equation*}
with $(x,t)=(x_1,\ldots, x_m,t)$, integers $a_{ij}$, and
$h_{ij}:K^m\to K$ $\cL$-definable functions which are $K$-analytic
on the base of $A_{i}$, $j=1,\ldots,r$. If $\lambda_i=0$, we use the
convention that $a_{ij}=0$. Moreover, given $\ell,n>0$, we can take the
$A_i$ such that moreover
 $$f_j(x,t)\cdot Q_{\ell,n}$$
for $(x,t)\in A_i$ only depends on $i$ and $j$ (and not on
$(x,t)$), and such that the restriction of $f_j$ to each $A_i$ is $K$-analytic.
 \end{theorem}

\begin{definition}\label{full}
If $f_j$ and the $A_i$ are as in Theorem \ref{thm:CellDecomp}, then
call $f_j$ \emph{prepared} on the cells $A_i$.
 If the base of $A_i$ is itself a cell on which the $h_{ij}(x)$ and the boundaries of $A_i$ are prepared, and so on $m$ times, then we call $A_i$ a \emph{full cell} and we call $f_j$ \emph{fully prepared} on the $A_i$. It is also clear what is meant by a full cell $A\subset K^m\times Y$ \emph{over} some $\cL$-definable set $Y$, in analogy to the notion of cells over $Y$ of Definition \ref{def::cell}.
\end{definition}

Clearly by induction on  dimension (that is, on $m$) one can use Theorem \ref{thm:CellDecomp} to get a partition into full cells on which the $f_j$ are fully prepared.

We formulate four basic corollaries of Theorem \ref{thm:CellDecomp}.
 The first one was originally
proven without using Theorem \ref{thm:CellDecomp} in
\cite{vdDriesSkolem} and \cite{DvdD}.

\begin{cor}[Definable Skolem functions]\label{skolem}
Let $X\subset K^n\times K^m$ be an $\cL$-definable set. Then there
exists an $\cL$-definable function $f:K^n\to K^m$ such that for each
$(x,y)\in X$  with $x\in K^n$ and $y\in K^m$ the point $(x,f(x))$
lies in $X$.
\end{cor}

\begin{cor}[Uniform boundedness]\label{finite}
Let $X\subset K^{n}\times K^m$ be $\cL$-definable, with $n,m\geq 0$. Then there exists $N>0$
such that for all $y\in K^m$ with $X_y:=\{x\in K^n\mid (x,y)\in X\}$ a finite set, one has
$$
\# X_y < N.
$$
Moreover, any discrete $\cL$-definable set $A\subset K^n$ is finite.
\end{cor}

\begin{cor}[Injectivity versus constancy]\label{inj-cons}
Let $Y$ and $X\subset K\times Y$ be $\cL$-definable sets and let
$F:X\to K$ be an $\cL$-definable function. Then there exists a
finite partition of $X$ into $\cL$-definable sets $X_i$ such that
for each $y\in Y$, the restriction of $F(\cdot,y): x \mapsto F(x,y)$
to
$$
X_{iy}:=\{x\in K\mid (x,y)\in X_i\}
$$
is either injective or constant, where this distinction only depends
on $i$ (and not on $y$).
\end{cor}
\begin{proof}

Let $\Gamma\subset K^2\times Y$ be the graph of $F$. Now let
$p:\Gamma\to K\times Y$ be the coordinate projection sending
$(x,F(x,y),y)$ to $(F(x,y),y)$ and let $p':\Gamma\to K\times Y$ be the projection
sending $(x,F(x,y),y)$ to $(x,y)$. Apply Theorem
\ref{thm:CellDecomp} to $\Gamma$ over $K\times Y$ according to $p$
(that is, the cells have a cell-like description in the
$x$-variable). For each $1$-cell $\Gamma_i$ in the partition of
$\Gamma$, the corresponding restriction of $F(\cdot,y)$ to
$p'(\Gamma_i)_y$ is clearly locally constant.
 For each $0$-cell $\Gamma_i$ in the partition of
$\Gamma$,  the corresponding restriction of $F(\cdot,y)$ to
$p'(\Gamma_i)_y$ is clearly injective. One completes the
proof by Corollaries \ref{finite} and \ref{skolem}.
\end{proof}

\begin{cor}[Cell criterion]\label{ortho2}
Let
$Y$ and
$X\subset K\times Y$ be $\cL$-definable and let $d:Y\to K$ be an $\cL$-definable function
and $n>0$. Suppose that for each $(t,y)\in X$ with $y\in Y$ there is a maximal ball $B_{t,y}$ containing $t$ such that $B_{t,y}\times \{y\}$ is contained in $X$. Suppose further that
$$
B_{t,y}=\{w\in K\mid  \ord(w-d(y)) = b_{t,y},\ \ac_{n}(w-d(y))=\xi_{t,y}\}
$$
for some $b_{t,y}$ and $\xi_{t,y}\not=0$. Then $X$ is a finite disjoint union of $p$-adic $1$-cells $A_i$ with center the restriction of $d$ to the base of $A_i$, and such that each ball $B_{t,y}$ appears as a ball above $y$ of one of the cells $A_i$.
\end{cor}
\begin{proof}
Since the image of  $\ac_{n}$ is finite we may suppose that $\xi_{t,y}$ is constant. Now the corollary follows
from Theorem \ref{thm:CellDecomp} and from Presburger cell decomposition results of \cite{pres} in a straightforward way.
\end{proof}

\begin{lem}\label{critII}
Let
$Y$ and
$X\subset K\times Y$ be $\cL$-definable.  Suppose that for each $(t,y)\in X$ with $y\in Y$ there is a maximal ball $B_{t,y}$ containing $t$ such that $B_{t,y}\times \{y\}$ is contained in $X$. Then $X$ is a finite disjoint union of $p$-adic $1$-cells $A_i$ such that each ball $B_{t,y}$ appears as a ball  above $y$ of one of the cells $A_i$.
\end{lem}
\begin{proof}
Note that each ball $B_{t,y}$ can be written as
\begin{equation}\label{eq:II}
B_{t,y}=\{z\in K\mid \ord (z-w)= a_{t,y},\ z-w\in Q_{1,1}\},
\end{equation}
for unique $a_{t,y}\in \ZZ$ and for (non unique) $w\in K$.
In the basic case that for each $y$ in $Y$ the set $X_y:=\{t\in K\mid (t,y)\in X\}$ is a ball one automatically has that $X_y=B_{t,y}$ for all $(t,y)\in X$. In this basic case define $W$ as
$$
W=\{(w,y)\in K\times Y\mid \mbox{ Equation (\ref{eq:II}) holds for $B_{t,y}$ and $w$} \}.
$$
Now use Corollary \ref{skolem} to find an $\cL$-definable function $d:Y\to K$ whose graph lies in $W$ and use Corollary \ref{ortho2} to finish this basic case.

In the general case partition $X$ into finitely many cells $X_i$ over $Y$ with center $c_i$, coset $\lambda_iQ_{m_i,n_i}$, and base $Y_i$ by using Theorem \ref{thm:CellDecomp}. Up to refining the partition $\{X_i\}$ of $X$, we may suppose that the following distinction only depends on $i$ when $y$ moves over $Y_i$: either $c_i(y)$ lies inside $B_{t,y}$ for some $t$, or, $c_i(y)$ lies outside $B_{t,y}$ for all $t$. Define $I_1$ and $I_2$ such that $i\in I_1$ if and only if $c_i(y)$ lies inside $B_{t,y}$ for some $t$, and $i\in I_2$ else.
Then, for $i\in I_1$ and $y\in Y_i$, let $B(i,y)$ be the ball $B_{t,y}$ containing $c_i(y)$ and define
$$
A_i:=\{(t,y)\mid y\in Y_i,\ t \in B(i,y)\}.
$$
Then the $A_i$ for $i\in I_1$ are as in the basic case and can thus be treated.
Put
$$
X' = X\setminus ( \bigcup_{i\in I_1} A_i ).
$$
It is enough to prove the statement of the lemma for $X'$.
Fix $(t,y)\in X'$ and choose $i$ (either in $I_1$ or in $I_2$) such that $B_{t,y}$ contains
at least one of the balls of $X_i$ above $y$. (Such $i$ must exist by the non archimedean property and the maximality of the occurring balls.)
Since by construction $c_i(y)$ lies outside $B_{t,y}$ for the fixed $(t,y)$,
there are $b_{t,y}\in\ZZ$, $m\leq m_i$, and $\lambda \in K^\times$
such that
$$
B_{t,y} = \{z\in K\mid \ord ( z - c_i(y) )= b_{t,y},\ z - c_i(y) \in \lambda Q_{m,n_i} \}.
$$
Since there are only finitely many $i$, $m$, and cosets of $Q_{m,n_i}$ in $K^\times$, we can finish by Corollary \ref{ortho2}.
\end{proof}

\begin{definition}[Jacobian property]\label{defjacprop}
Let $F:B_1\to B_2$ be an $\cL$-definable function with
$B_1,B_2\subset K$. Say that $F$
\textit{has the Jacobian property} if the following
conditions a) up to d) hold
\begin{itemize}

\item[a)] $F$ is a bijection $B_1\to B_2$ and $B_1$ and $B_2$ are balls;

\item[b)] $F$ is $C^1$ on $B_1$; write ${\rm Jac} F$ for $\partial
F/\partial x:B_1\to K$;

\item[c)] $\ord ({\rm Jac} F)$ is constant (and finite) on $B_1$;

\item[d)] for all $x,y\in B_1$ with $x\not=y$, one has
$$
\ord ({\rm Jac} F)+\ord(x-y)=\ord(F(x)-F(y)).
$$
\end{itemize}
\end{definition}

\begin{prop}[Jacobian property for definable functions \cite{CLip}, Section 6]\label{jacprop}
Let $Y$ and $X\subset K\times Y$ be $\cL$-definable sets, let
$F:X\to K$ be an $\cL$-definable function.
Suppose that for each $y\in Y$, the function $F(\cdot,y):t\mapsto F(t,y)$ is injective.
 Then there exists a finite partition of $X$ into $p$-adic
cells $A_i$ over $Y$ such that for each $i$, each $y\in Y$ and each
ball $B$ of $A_i$ above $y$, there is a ball
$B^* \subset K$
such that the map
$$
F_B:B\to B^*:t\mapsto F(t,y)
$$
is well defined and has the Jacobian property.
 \end{prop}

The following proposition is new and relies on Proposition \ref{jacprop}.

\begin{prop}[Compatible cell decompositions under a definable function $F$]\label{celcel}
Let $X$, $Y$, and $F$ be as in Proposition \ref{jacprop}, where in particular  $F(\cdot,y)$ is injective for each $y\in Y$. Define
$F_Y$ as the ($\cL$-definable) function $F_Y:X\to K\times Y:(t,y)\mapsto
(F(t,y),y)$.
 For $A_i$ as in Proposition \ref{jacprop}, write $A'_i$ for the set
$F_Y(A_i)$.
Then we can choose the partition of $X$ into cells $A_i$ over $Y$ as
in Proposition \ref{jacprop} such that moreover each $A'_i$ is a $p$-adic
cell over $Y$, and such that for each $y\in Y$ and each ball $B$ of
$A_i$ above $y$, $B^*$ is a ball of $A'_i$ above $y$, where $B^*$ is
as in Proposition \ref{jacprop}. Hence, for any $y\in Y$, there
is a correspondence between the balls of $A_i$ above $y$ and the
balls of $A'_i$ above $y$.
\end{prop}

\begin{proof}
Partition $X$ into cells $X_i$ as in Proposition \ref{jacprop}. Up to this finite partition we may suppose that $X$
equals $X_1$ which we may suppose is a $1$-cell.   Write $X_1'$ for $F_Y(X_1)$.
Partition $X_1'$ into cells $X_{1i}'$ with center $d_i$, coset $\lambda_i Q_{m_i,n_i}$, and base $Y_i$ by using Theorem \ref{thm:CellDecomp}.
 For a ball $B=B_{t,y}$ of $A_1$ above $y$ containing $t$ write $B^*_{t,y}$ for the corresponding ball $B^*$, as given by Proposition \ref{jacprop}.
 Up to refining the partition $\{X_{1i}'\}$ of $X_1'$, we may suppose that the following distinction only depends on $i$ when $y$ moves over $Y_i$: either $d_i(y)$ lies inside $B^*_{t,y}$ for some $t$, or, $d_i(y)$ lies outside $B^*_{t,y}$ for all $t$.
Define $I_1$ and $I_2$ such that $i\in I_1$ if and only if $d_i(y)$ lies inside $B^*_{t,y}$ for some $t$, and $i\in I_2$ else.
Then, for $i\in I_1$ and $y$ in $Y_i$, let $B^*(i,y)$ be the ball $B^*_{t,y}$ containing $d_i(y)$. Let $B(i,y)$ be the ball $B_{t,y}$ where $t$ is such that $B^*_{t,y}=B^*(i,y)$. By construction,
\begin{equation}\label{FBiy}
F(B(i,y)\times \{y\}) = B^*(i,y).
\end{equation}
For $i\in I_1$, define
$$
A_i:=\{(t,y)\mid y\in Y_i,\ t \in B(i,y)\}
$$
and put $A_i':=F_Y(A_i)$. Apply Lemma \ref{critII} to $A_i$ and $A_i'$ for each $i\in I_1$.
Then these $A_i$ and $A_i'$ are as required by (\ref{FBiy}) and it is thus sufficient to prove the proposition for the restriction of $F$ to
$$
X':=X\setminus ( \bigcup_{i\in I_1} A_i ).
$$
 Fix $(t,y)\in X'$ and choose $i$ (either in $I_1$ or in $I_2$) such that $B^*_{t,y}$ contains
at least one of the balls of $X'_i$ above $y$. (Such $i$ must exist by the maximality of the occurring balls.)
Since by construction $d_i(y)$ lies outside $B^*_{t,y}$, and since $B^*_{t,y}$ contains
a ball of $X'_i$ above $y$,
there are $b_{t,y}\in\ZZ$, $m\leq m_i$, and $\lambda \in K^\times$
such that
$$
B^*_{t,y} = \{z\in K\mid \ord ( z - d_i(y) )= b_{t,y},\ z - d_i(y) \in \lambda Q_{m,n_i} \}.
$$
Since there are only finitely many cosets of $Q_{m,n_i}$ in $K^\times$, the proposition follows from Theorem \ref{thm:CellDecomp} and from the Presburger cell decomposition results of \cite{pres} in a straightforward way.
\end{proof}

Further we give a corollary of Proposition \ref{jacprop} that we will not use further on.

\begin{cor}\label{Lip-1}
Let
$Y$ be an $\cL$-definable set
$F:X\subset K\times Y\to K$ be an $\cL$-definable function such that $F(\cdot,y)$ is injective for each $y\in Y$. Then there exists a finite partition of $X$ into $\cL$-definable pieces $X_i$
such that, for each $i$ and each $y\in Y$, the restriction of $F(\cdot , y)$ to $X_{iy}:=\{t\in K\mid (t,y)\in X_i\}$ or its inverse function is locally $1$-Lipschitz.
\end{cor}
\begin{proof}
Apply Proposition \ref{jacprop} to $F$, yielding a partition of $X$ into cells $A_i$ over $Y$. Now partition each $A_i$ into pieces according to the condition that $|\partial F/\partial t|$ is $\leq 1$, resp.~$>1$ on the piece. On the pieces where $|\partial F/\partial t|$ is $\leq 1$ we are done by the Jacobian property which holds by construction. On a piece where $|\partial F/\partial t|$ is $> 1$, the inverse of $F(\cdot , y)$ is locally $1$-Lipschitz by the chain rule for differentiation and the Jacobian property which holds by construction.
\end{proof}

Note that the different possibilities for the (non exclusive) disjunctions in Corollary \ref{Lip-1} can be supposed to depend only on $i$ (and not on $y$) by taking the parts $X_i$ small enough. Indeed, the occurring conditions as local $1$-Lipschitz continuity, injectivity, and so on, are $\cL$-definable in $y\in Y$.

\section{Proofs of the main results}

Theorem \ref{rellip} and Proposition \ref{Lipcenter}
are proved using a joint induction on $m$.

\begin{proof}[Proof of Theorem \ref{rellip} for $m=1$]
We are given $\varepsilon>0$, $Y$ an $\cL$-definable set, and
$f:X\subset K\times Y\to K$ an $\cL$-definable function such that
for each $y\in Y$ the function $f(\cdot,y):x\mapsto f(x,y)$ is
locally $\varepsilon$-Lipschitz continuous on its natural domain
$X_y:=\{x\in K\mid (x,y)\in X\}$.
 Using Corollary \ref{inj-cons}, we may suppose that $f(\cdot,y)$ is
injective for each $y$.
 Use Proposition
\ref{celcel} to partition $X$ into finitely many $p$-adic cells
$X_i$ over $Y$ with center $c_i$.
 By working piecewise we may suppose that $X=X_1$ and that $X_1$ is a $1$-cell over $Y$.
 By the Jacobian property $f(\cdot,y)$ is $C^1$ and by local $\varepsilon$-Lipschitz
continuity,
\begin{equation}\label{jace} |\partial f(x,y)/\partial x
|\leq \varepsilon
\end{equation}
for all $(x,y)\in X$.
 By the above application of Proposition \ref{celcel}, the set
$$
X':=f_Y(X),
$$
with $f_Y:X\to K\times Y:(x,y)\mapsto (f(x,y),y)$, is a $p$-adic
$1$-cell with some center $d_1$.
 Since a function $g:A\subset K\to K$ is $C$-Lipschitz continuous if
and only if $A\to K:x\mapsto g(x+a)+b$ is $C$-Lipschitz continuous
for any constants $a,b\in K$, we may thus suppose, after translating, that $c_1$
and $d_1$ are identically zero.

Now fix $y\in Y$. Take $(x_1,y)$ and $(x_2,y)$ in $X$. If $x_1$ and
$x_2$ both lie in the ball $B_{x_1,y}$, then
\begin{equation}\label{m1jac}
|( \partial f(x_1,y)/\partial x)  \cdot (x_1-x_2)|= | f(x_1,y) -
f(x_2,y) |
\end{equation}
 by the Jacobian property. By (\ref{jace}) we are done and can take any $C\geq \varepsilon$.

Next suppose that $B_{x_1,y}$ and $B_{x_2,y}$  are two different
balls. By our assumption that $c$ and $d$ are identically zero, we can
write
$$
 B_{x_i,y} = \{x\in K\mid \ord (x)=a_{x_i,y},\ \ac_m(x)=\ac_m\lambda \}
$$
and likewise for their images under $f(\cdot,y)$,
$$
 B^*_{x_i,y} = \{z\in K \mid \ord (z)=b_{x_i,y},\ \ac_{m'}(z)=\ac_{m'}\mu
 \}.
$$
 From these descriptions we get the inequalities:
$$
\ord(f(x_1,y)- f(x_2,y)) = \min_{i=1,2} (b_{x_i,y})
$$
and
$$
\min_{i=1,2} (a_{x_i,y} ) =   \ord(x_1-x_2)
$$
 On the other hand by the very Jacobian property d) one has
$$
m +  \ord (\partial f(x_i,y)/\partial x)  + a_{x_i,y}  =  m' +
b_{x_i,y}
$$
for $i=1,2$. Hence, putting this together with (\ref{jace}),
$$
 |f(x_1,y)- f(x_2,y)| = \max_{i=1,2} p^{-b_{x_i,y}}
 \leq   \varepsilon  p^{m'-m}   \max_{i=1,2} p^{-a_{x_i,y}}  =  \varepsilon p^{m'-m} |x_1-x_2|
$$
and thus one can take any $C\geq \varepsilon\cdot \max(1,
p^{m'-m})$.
\end{proof}

\begin{remark}\label{chain}
The chain rule for
differentiation yields the following statement for a $C^1$ function $f:X\subset K\times Y\to K$ on an open set $X$ in the variables
$(t,y)$, where $y=(y_1,\ldots,y_m)$ runs over an open $Y\subset K^m$ and $t$ over $K$. Suppose that 
the function $f(\cdot,y):t\mapsto f(t,y)$ is injective and has $C^1$ inverse for each $y\in Y$. Define $Z$ as $\{(f(t,y),y)\in
K\times Y \mid  (t,y)\in X\}$ and define the function
$g:Z\to K$ by $(z,y)\mapsto t$ for the unique $t$ with $f(t,y)=z$. Then one has for
each $i=1,\ldots,m$
$$
\frac{\partial g(z,y)}{\partial y_i} = - \frac{\partial
f(t,y)}{\partial y_i} \cdot \Bigl(\frac{\partial f(t,y)}{\partial
t}\Bigr)^{-1}
$$
where $z=f(t,y)$.
\end{remark}

\begin{proof}[Proof of Proposition \ref{Lipcenter} for $m$ using Theorem
\ref{rellip} for $m-1$.]

We will proceed by induction on $m$. For $m=1$ the statement of Proposition \ref{Lipcenter} is trivial and hence we may suppose
that $m>1$.
 Up to a finite partition of $X$, we may assume that $X$ itself is
either a $1$-cell or a $0$-cell over $K^{m-1}\times Y$ along some
coordinate projection $p:K^m\times Y\to K^{m-1}\times Y$, say with
center $c$.

\par
First suppose we are in the basic case that, for $y$ such that $X_y$ is nonempty, the set $p(X)_y$ is not open in $K^{m-1}$.
By the induction hypotheses, we may suppose  that $p(X)$ is a $p$-adic cell (over $K^{m-2}\times Y$) with center $c_{m-1}$ such that $c_{m-1}(\cdot,y))$ is $C$-Lipschitz in the (relevant) variables $x_1,\ldots,x_{m-2}$ for each fixed value of $y\in Y$,
and so on for in total $m-1$ subsequent coordinate projections, up to the projection to $Y$. Then, after the (triangular) bi-Lipschitz continuous transformation where we replace $x_{m-1}$ by $x_{m-1} - c_{m-1}$ and so on $m-1$ times, we may suppose that the center of the cell $p(X)$ is identically zero, and so on $m-1$ times up to the projection to $Y$. If we still
 use
 the name $p(X)$ for the so-obtained transformed set, there must be a coordinate $x_i$, for some $i=1,\ldots,m-1$, which is identically zero on $p(X)$, and thus we can finish by Proposition \ref{Lipcenter} for $m-1$.

\par

Hence, we can place ourselves in the more interesting case that the $p(X)_y$ are open in $K^{m-1}$ for all $y\in Y$ and of course we may then suppose that moreover $c(\cdot,y)$ is $C^1$ on $p(X)_y$ for all $y\in Y$.
 After reordering the variables $x_1,\ldots,x_{m-1}$ and after finitely partitioning $p(X)$, we may suppose that $|\partial c/\partial x_{m-1}|$ is maximal among the $|\partial c/\partial x_{i}| $ on the whole of $p(X)$ for $i=1,\ldots,m-1$. If $|\partial c/\partial x_{m-1}|\leq 1$ on the whole of $p(X)$, then we are done by Proposition \ref{localL} and Theorem \ref{rellip} for $m-1$, up to a finite partition of $p(X)$. Hence, we may further assume, up to a finite partition of $p(X)$, that
$1< |\partial c/\partial x_{m-1}|$ on the whole of $p(X)$.
 Using Corollary \ref{inj-cons} and up to a further finite partition of $p(X)$, we may furthermore suppose that
$c(x_1,\ldots,x_{m-2},\cdot,y)$ is injective for each
$(x_1,\ldots,x_{m-2},y)$.
Now partition $p(X)$ again, as follows. Use Proposition \ref{celcel} for the map $c(x_1,\ldots,x_{m-2},\cdot,y)$ to partition $p(X)$ into
finitely many $p$-adic cells $A_i$ over $K^{m-2}\times Y$, along the
projection $(x_1,\ldots,x_{m-1},y)\mapsto (x_1,\ldots,x_{m-2},y)$.
Up to such a finite
partition of $p(X)$, we may suppose that $p(X)=A_1$ and that $A_1$ has center $c_1$.

 \par
First we treat the more simple case that $X$ is a $0$-cell over $p(X)$.
 In this case we simply invert the role of $x_m$ and $x_{m-1}$ in the
build-up of the cell $X$ as follows. Write
$d(x_1,\ldots,x_{m-2},\cdot,y)$ for the inverse function
$(c(x_1,\ldots,x_{m-2},\cdot,y))^{-1}$ of
$c(x_1,\ldots,x_{m-2},\cdot,y)$.
 Then $X$ is also a $0$-cell over $K^{m-1}\times Y$ along the projection
$p'$ sending $(x_1,\ldots,x_{m},y)$ to $(x_1,\ldots,x_{m-2},x_m,y)$
with center $d$. Since $d$ is constructed as an inverse function and if we recall the differentiation rule for inverse functions and Remark \ref{chain}, it is clear that all partial derivatives $\partial d/\partial x_i$ for $i=1,2\ldots,m-2,m$ are bounded in norm. We are done by Proposition \ref{localL} and Theorem \ref{rellip} for $m-1$, up to a finite partition of $p'(X)$.

\par
Finally,  we treat the most interesting case that $X$ is a $1$-cell.
 For $1\leq i\leq m$ and $a\in K^m$ write $\hat a_i$ for
$(a_1,\ldots,a_i)$.
 Fix $(a,y)$ in $X$. Let $B_{a,y}$ be the unique ball (of the cell
$X$) above $(\hat a_{m-1},y)$ which contains $a_m$. Further, let $B^0_{a,y}$ be the unique ball of the cell $p(X)$ that contains
$a_{m-1}$ and lies above $(\hat a_{m-2},y)$.
 By the previous application of Proposition \ref{celcel} for the map $c(x_1,\ldots,x_{m-2},\cdot,y)$, the image
of $B^0_{a,y}$ under $c(\hat a_{m-2},\cdot,y)$ is a ball
$B^{0*}_{a,y}$ and one has moreover descriptions, uniformly in
$(a,y)$ in $X$,
$$
B_{a,y}=\{x_m\mid \ord (x_m-c(\hat a_{m-1},y ) ) = b_{a,y}, \
\ac_{n} (x_m-c(\hat a_{m-1},y ) ) = \ac_{n} \lambda \}
$$
$$
B^0_{a,y} =\{x_{m-1}\mid \ord (x_{m-1}-c_1(\hat a_{m-2},y )
) = b^0_{a,y}, \ \ac_{n'} (x_{m-1}-c_1(\hat a_{m-2},y )
) = \ac_{n'} \lambda' \}
$$
and
$$
B^{0*}_{a,y} = \{ z\mid \ord (z-e(\hat a_{m-2},y ) ) =
b^{0*}_{a,y}, \ \ac_{n''} (z-e(\hat a_{m-2},y ) ) =
\ac_{n''} \lambda'' \},
$$
for some nonzero constants $\lambda$, $\lambda'$, $\lambda''$, $n$, $n'$,
$n''$ coming from the descriptions of the cells, where $e$ is the $\cL$-definable function as given by the previous
application of Proposition \ref{celcel}, and where $b_{a,y}$ only
depends on the ball $B_{a,y}$, and similarly $b^0_{a,y}$ only depends on $B^0_{a,y}$ and
$b^{0*}_{a,y}$ only on $B^{0*}_{a,y}$.

We will compare sizes of balls, where we call a ball $B_1$ strictly bigger in size than a ball $B_2$ if a translate of $B_2$ is strictly contained in $B_1$ and we say that $B_1$ and $B_2$ are equal in size of a translate of $B_1$ equals $B_2$.
By partitioning $p(X)$ further we may assume that we are in one of the following two cases.\\

\textbf{Case 1. The ball $B_{a,y}$ is bigger or equal in size than  $B^{0*}_{a,y}$ for all $(a,y)$ in $X$.} \\

Case 1 is equivalent
to
 $b_{a,y}^{0*} \geq n - n'' + b_{a,y}$ for all $(a,y)$ in $X$.
Also, $B_{a,y}$ does not depend on $a_{m-1}$ when $a_{m-1}$
runs over $B^0_{a,y}$, but the center $c$ itself may of course depend nontrivially on $x_{m-1}$. We will replace the center $c$ by another center which depends trivially on $x_{m-1}$, as follows.
   By construction $c(\hat a_{m-1},y)$ lies in $B^{0*}_{a,y}$ and the set $B^{0*}_{a,y}$ is described above.
By this description, $e$ is a kind of approximation of $c$ and is thus a candidate to become the new center instead of $c$, which we show indeed to work as follows.
  Let $\ell$ be $\max(n,n'')$. Partition $X$ into finitely many parts where $\ac_{\ell}( x_m - e(\hat a_{m-2},y ) )$ is constant. Next, apply Corollary \ref{ortho2} to each such part to obtain a partition of $X$ into finitely many cells $A_{i}$ with center the restriction of $e$ to the base of $A_{i}$ (Corollary \ref{ortho2} can be applied because $\ell$ is well chosen).
  Hence, up to this finite partition we can suppose that
the center of the $p$-adic cell $X$ does not depend on the variable $x_{m-1}$. With this new situation, we can go back in the proof and start reordering the variables $x_1,\ldots,x_{m-1}$ and finitely partitioning $p(X)$ such that $|\partial c/\partial x_{m-1}|$ is again maximal among the $|\partial c/\partial x_{i}| $ on the whole of $p(X)$ for $i=1,\ldots,m-1$, as we did above. After finitely many recursions, we will not fall into Case 1 anymore. Indeed, if the $|\partial c/\partial x_{i}| $ are $\leq 1$ then we are in a case treated above and so on.  \\

\textbf{Case 2. The ball $B^{0*}_{a,y}$  is strictly bigger in size than $B_{a,y}$ for all $(a,y)$ in $X$.}\\

Case 2 is equivalent
to
$ n + b_{a,y} >  n'' + b_{a,y}^{0*}$ for all $(a,y)$ in $X$ and implies
\begin{equation}\label{BinB}
B_{a,y} \subset B^{0*}_{a,y} \quad \mbox{ for all $(a,y)\in X$.}
\end{equation}
 By construction, we can consider the inverse function of $c(\hat a_{m-2},\cdot,y)$.
 Write
$d(\hat a_{m-2},\cdot,y)$ for the inverse function of $c(\hat a_{m-2},\cdot,y)$. Then the domain of $d(\hat a_{m-2},\cdot,y)$ contains in particular the ball $B_{a,y}$ by (\ref{BinB}) and hence, we can apply $d(\hat a_{m-2},\cdot,y)$ to $x_m$ for any point $(\hat a_{m-2},x_{m-1},x_m,y)$ in $X$.
 Partition $X$ into finitely many parts where $\ac_{\ell'}(x_{m-1} - d(\hat a_{m-2},x_m,y ) )$ is constant for some sufficiently large $\ell'$. Since $\ell'$ is sufficiently large and by the Jacobian property which holds by the previous application of Proposition \ref{celcel}, we can apply Corollary \ref{ortho2} to each such part to obtain a partition of $X$ into finitely many $p$-adic cells $A_{i}$ with center the restriction of $d$ to the base of $A_{i}$.
 Up to this partition, we may suppose that $X$ is a $p$-adic cell over $K^{m-1}\times Y$
along the projection $p'$ sending $(x,y)$ to $(\hat x_{m-2},x_m,y)$ with center $d$.
Since $d$ is constructed as an inverse function and if we recall the differentiation rule for inverse functions and Remark \ref{chain}, it is clear that all partial derivatives $\partial d/\partial x_i$ for $i=1,2,\ldots,m-2,m$ are bounded in norm. Hence we can finish by Proposition \ref{localL} and Theorem \ref{rellip} for $m-1$, up to a finite partition of $p'(X)$.
\end{proof}

\begin{proof}[Proof of Theorem \ref{rellip} for $m>1$, using Proposition \ref{Lipcenter} for $m$.]
We proceed by induction on $m$, where the case $m=1$ is proven above in this section.
 We are given $\varepsilon>0$, an $\cL$-definable set $Y$ and an
$\cL$-definable function $f:X\subset K^m\times Y\to K$ such
that for each $y\in Y$ the function $f(\cdot,y):x\mapsto f(x,y)$ is
locally $\varepsilon$-Lipschitz continuous.
Merely to ensure later on that partial derivatives are well defined when they will appear, we may now already suppose that the  $X_y$ are $K$-analytic manifolds on which $f(\cdot,y)$ is $K$-analytic for each $y\in Y$, but we will not necessarily maintain this property throughout the proof for the parts of upcoming partitions.
 Use the notation $\hat x$ for $(x_1,\ldots,x_{m-1})$ and similarly for tuples $a=(a_1,\ldots,a_{m})$ in $K^m$ for which $\hat a=(a_1,\ldots,a_{m-1})$.
 By Theorem \ref{rellip} for $m-1$ and up to a finite partition of
$X$, we may suppose that for each $(a,y)=(a_1,\ldots,a_m,y)$
in $X$ each of the functions
\begin{equation}\label{PL}
 f(\cdot,a_i,\cdot,y):
(x_1,\ldots,x_{i-1},x_{i+1},\ldots,x_m)\mapsto  f(x_1,\ldots,x_{i-1},a_i,x_{i+1},\ldots,x_m,y)
\end{equation}
is $C$-Lipschitz continuous.
  By Theorem \ref{thm:CellDecomp} we may suppose that $X$ is a cell over $K^{m-1}\times Y$, say, with center $c$, along a coordinate projection $p:X\to K^{m-1}\times Y$.
 By Proposition \ref{Lipcenter} for $m$, we may suppose that $c(\cdot,y)$ is $C$-Lipschitz continuous in $\hat x$ for each $y\in Y$. Up to a finite partition of $p(X)$, we may
 suppose
 by Theorem \ref{rellip} for $m-1$  that, for each fixed values of $a_m$ and $y$, the function $\hat x\mapsto f(\hat x, a_m-c (\hat x,y),y)$ is $C$-Lipschitz continuous.
 Hence, if we perform the bi-Lipschitz transformation which replaces $x_m$ by $x_m-c (\hat x,y)$ but which preserves the other coordinates, then we see that we may suppose that:
\begin{itemize}
\item[($*$)] $X$ is a cell over $K^{m-1}\times Y$ whose center is identically zero, the function $f(\cdot,a_m,y)$ is $C$-Lipschitz continuous in $\hat x$ for each $a_m\in K$ and $y\in Y$, and for all fixed $\hat x, y$, the function $f(\hat x,\cdot, y)$ is also $C$-Lipschitz continuous in $x_m$.
\end{itemize}
 In the simple case that $X$ is a $0$-cell over $K^{m-1}\times Y$, we can finish by Theorem \ref{rellip} for $m-1$, since the $x_m$-coordinate is identically zero on $X$ and can be neglected.
 Next suppose that $X$ is a $1$-cell.
 Let $(a,y)$ and $(a',y)$ be given in $X$.
 If $|a_m|= |a_m'|$, then the point $(\hat a,a'_m,y)$ also lies in $X$ by the definition of cells and since the center of $X$ is zero, and hence we can jump inside $X$ from $(a,y)$ to $(\hat a,a'_m,y)$ and finally jump further to $(a',y)$. Calculating the images under $f$ and controlling the distances between these points  yields:
$$
|f(a,y)-f(a',y)|
$$
 $$
= |f(a,y)-f(\hat a,a'_m,y) +
f(\hat a,a'_m,y)-f(a',y)|
$$

$$
\leq  \max (|f(a,y)- f(\hat a,a'_m,y) |,\
|f(\hat a,a'_m,y)-f(a',y)|)
$$
 $$
\leq \max ( C|a_m - a_m' |,\  C|(\hat a,a'_m) - a'|)
$$
$$
= C|a - a'|,
$$
where the first inequality follows from the non archimedean property, the second from property ($*$), and the last equality follows from the definition of the norm on $K^m$.
 Let us now suppose that $|a_m|\not = |a_m'|$, say, $|a'_m|<|a_m|$.
  First suppose that $|x_m|$ has no lower bound in $X$, that is, for each $(\hat x,x_m,y)$ in $X$ there exists $x_m'\in K$ arbitrarily close to $0$ such that
$(\hat x,x_m',y)$ still lies in $X$. Then again the point $(\hat a,a'_m,y)$ lies in $X$ by the definition of cells and since the center of $X$ is zero. Hence we can make the same jumps as in the previous case and the same computation will hold for the same reasons.
 In the other case we may suppose that for each $(b,y)$ in $X$ there
is a minimal value $e(\hat b,y)>0$ among the values $|x_m|$
for all $x_m$ with  $(\hat b,x_m,y)$ in $X$. Let $\alpha(\hat
b,y)\in K$ be such that $|\alpha(\hat b,y)| =
e(\hat b,y)$. By Corollary \ref{skolem} we may suppose that $\alpha$
is an $\cL$-definable function in $(\hat b,y)$ whose graph lies in $X$.

\par

Up to a finite partition we may
suppose,
by Proposition \ref{Lipcenter} for $m-1$,
that $p(X)$ is a full cell over $Y$ whose centers are all $C$-Lipschitz continuous. If for $y\in Y$, $X_y$ is not open in $K^m$ if nonempty, then we can can, after a natural triangular transformation, force one of the coordinates $x_i$ for $i=1,\ldots,m$ to be $0$ on $X$, as in the proof of Proposition \ref{Lipcenter}. Hence in this case, we are again done by Theorem \ref{rellip} for $m-1$. Thus, we may suppose that $X_y$ is open in $K^m$ for each $y$.

\par

Up to a finite partition of $p(X)$ and by Proposition \ref{localL}, we may suppose that we are in one of the two following cases. \\

\textbf{Case 1. The function $\alpha$ is locally $C$-Lipschitz continuous on $p(X)_y$ for each $y\in Y$.}\\

Case 1 implies that the functions $\hat x\mapsto f(\hat x,\alpha(\hat x),y)$ are also locally $C'$-Lipschitz continuous for some $C'$. By Theorem \ref{rellip} for $m-1$ and up to a finite partition of $p(X)$ we may suppose
 that, for all $y$,
\begin{equation}\label{eqL}
\hat x\mapsto f(\hat x,\alpha(\hat x),y)
\end{equation}
is $C$-Lipschitz continuous on the whole of $p(X)$, possibly by replacing $C$ by some bigger constant as
 allowed.
Omitting  $y$ out of the notation, we will jump from $a$ to $(\hat a,\alpha(\hat a))$, jump further to  $(\hat a',\alpha(\hat
a'))$, and finally to $a'$, where $\hat
a'$ of course stands for $(a_1',\ldots,a_{m-1}')$.  We
compute, still omitting  $y$ out of the notation,
 $$
|f(a)-f(a')|
$$
$$
= |f(a)-f(\hat a,\alpha(\hat a)) + f(\hat
a,\alpha(\hat a))-
f(\hat a',\alpha(\hat a'))+ f(\hat a',\alpha(\hat
a')) - f(a') |
$$
$$
\leq C \max\Big(   |a_m - \alpha(\hat a) | ,
|(\hat a,\alpha(\hat a)) -
    (\hat a',\alpha(\hat a')) |,
  | \alpha(\hat a') - a'_m | \Big)
$$
$$
 = C|a - a'|,
$$
where the inequality holds by the $C$-Lipschitz continuity of the function  (\ref{eqL}) and by property ($*$), and the last equality holds by properties of the non archimedean norm on $K^m$ and the facts that $|a_m|\not = |a_m'|$ and $|\alpha(\hat a)|\leq |a_m|$ and $|\alpha(\hat a)|\leq |a_m'|$.\\

\textbf{Case 2. $|\frac{\partial \alpha } { \partial x_i }| > 1 $ on $p(X)$ for some $i<m$.}\\

We may further suppose that for a certain $j$,  $|\partial \alpha /\partial x_j|$ is maximal
among the $|\partial \alpha / \partial x_i|$ on the whole of $p(X)$.
For notational simplicity, suppose that $j=m-1$ (the case $j<m-1$ is only notationally different). Hence, by Corollary \ref{inj-cons} and Proposition
\ref{jacprop} we may suppose that
$$
\alpha(\hat x_{m-2},\cdot,y)
$$
is injective and $C^1$ with $C^1$ inverse for each $(x,y)\in X$.
 Let $\beta(\hat x_{m-2},\cdot,y)$ be the inverse of $ \alpha(\hat
x_{m-2},\cdot,y)$. We will make exactly the same jumps as in case 1, and establish exactly the same series of inequalities,
but these inequalities will hold for different reasons as in case 1.
 Write $\hat x_{m-2}$ for $(x_1,\ldots,x_{m-2})$.
By
Remark \ref{chain}, by the differentiation rule for inverse functions, and by the above supposition that $f(\cdot,y)$ is $K$-analytic on $X_y$ for each $y$ (which still may be supposed to hold here), the function
$$
F:(\hat x_{m-2},x_m,y)\mapsto f(\hat x_{m-2},\beta(\hat
x_{m-2},x_{m},y),x_m ,y)
$$
has bounded partial derivatives along $x_i$ for $i=1,\ldots,m-2,m$, on its natural domain $X'\subset
K^{m-1}\times Y$ (where $X'_y$ is also open for each $y$, hence the
partial derivatives are well defined on $X'_y$). By Proposition \ref{localL} and Theorem \ref{rellip} for $m-1$, after
a finite partition of $p(X)$ (and hence of $X'$), we may suppose
that $F(\cdot,y)$ is $C$-Lipschitz continuous on $X'_y$ for each
$y$. Write $d$ for $\alpha(\hat a_{m-1},y)$ and $d'$ for
$\alpha(\hat a'_{m-1},y)$. Now we can compute, since $(\hat
a_{m-1},d,y),y)=(\hat a_{m-2},\beta(\hat a_{m-2},d,y),d,y)$ and
again omitting  $y$ from the notation,
 $$
|f(a)-f(a')|
$$
$$
= |f(a)-f(\hat a_{m-1},d) + f(\hat a_{m-2},\beta(\hat a_{m-2},d),d)-
$$
$$
f(\hat a_{m-2}',\beta(\hat a'_{m-2},d'),d')+ f(\hat
a_{m-2}',\beta(\hat a'_{m-2},d'),d') - f(a') |
$$
$$
\leq C \max\Big(   |a_m - d | ,
|(\hat a_{m-2},\beta(\hat a_{m-2},d),d) -
    (\hat a_{m-2}',\beta(\hat a'_{m-2},d'),d') |,
  | d' - a'_m | \Big)
$$
$$
 = C|a - a'|,
$$
where we have used that $F(\cdot,y)$ is $C$-Lipschitz continuous for
each $y$ (instead of the Lipschitz continuity of (\ref{eqL}) used in case 1) but further similar reasons as in case 1. Indeed, the above equations and inegalities have exactly the same meaning as in case 1, they are only written differently to make it apparent that the Lipschitz continuity of $F(\cdot,y)$ can be used.
\end{proof}

\subsection*{Elementary equivalent fields}

We conclude with an analogue of Theorem \ref{rellip} for $p$-adically closed fields. 
Let $K_1$ be a field which is elementary equivalent to $K$ in the language $\cL$. Then $K_1$ is a valued field and we write $\cO_{K_1}$ for its valuation ring. One uses the norm notation $|\cdot|$ for the natural map from $K$ to the ordered multiplicative semi-group $\Gamma_1:=K_1/(\cO_{K_1}^\times)$.
Using this norm, for rational $C>0$, one can take the obvious definition for Lipschitz continuity with constant $C$. More generally, for nonzero $C$ in the divisible hull of $\Gamma_1$, there is a natural notion of Lipschitz continuity with constant $C$.
  Note that in the following result we can take $C$ to be a rational number, which is stronger than allowing nonzero $C$ from $\Gamma_1$.

\begin{prop}\label{ext}
Let a rational number $\varepsilon>0$ be given. Let
$f:X\subset K_1^m\times Y\to K_1$ be an $\cL(K_1)$-definable function (that is, $\cL$-definable with parameters from the field $K_1$) such
that for each $y\in Y$ the function $f(\cdot,y):x\mapsto f(x,y)$ is
locally $\varepsilon$-Lipschitz continuous on $X_y$, where also $Y$ is $\cL(K_1)$-definable. Then there exist a rational number $C>0$
and a finite partition of $X$ into $\cL(K_1)$-definable parts $A_i$ such that
for each $y\in Y$ and $i$ the restriction of $f(\cdot,y)$ to
$A_{iy}$ is (globally) $C$-Lipschitz continuous.
\end{prop}
\begin{proof}
(The proof uses a standard technique for using
 results like Theorem \ref{rellip} which hold for definable families.)
In both the $\cL(K_1)$-formulas $\varphi_X$ and $\varphi_f$ describing $X$ and $f$ there appear only finitely many parameters from $K_1$, say $r_1,\ldots,r_s\in K_1$. Replace these parameters $r_1,\ldots,r_s$ by new variables, say $z_1,\ldots,z_s$ (that is, the $z_i$ do not yet play a role in the formulas $\varphi_X$ and $\varphi_f$). Let the tuple of variables $z=(z_1,\ldots,z_s)$ run over $K_1^s$. The obvious variants of $\varphi_X$ and $\varphi_f$, with the $r_i$ replaced by the $z_i$, are of course interpretable in the standard $p$-adic field $K$ itself. Now it is an exercise to construct $\cL$-definable families $\tilde X_{\tilde y}\subset K^m$  of sets and functions $\tilde f_{\tilde y}:\tilde X_{\tilde y}\to K$, for some suitable parameter $\tilde y$ (containing in particular the $z$-tuple), which fall under the conditions and thus the conclusion of Theorem \ref{rellip} in such a way that Proposition \ref{ext} follows when one fills in the values $r_i$ back in for $z_i$.
\end{proof}

\subsection*{Acknowledgment}
The authors would like to thank J.~Denef and L.~van den Dries for
interesting discussions on the subject of Lipschitz continuity
during the preparation of this paper.

\bibliographystyle{amsplain}

\begin{thebibliography}{SGA}



\bibitem{pres}
R. Cluckers, \textit{Presburger sets and $p$-minimal fields}, J.
Symbolic Logic \textbf{68} (2003),  153--162.


\bibitem{Ccell}
R. Cluckers, \textit{Analytic $p$-adic cell decomposition and
integrals}, Trans. Amer. Math. Soc. \textbf{356} (2004), 1489--1499.



\bibitem{CCL}
R. Cluckers, G. Comte, F. Loeser, \textit{Local metric properties of $p$-adic definable sets}.

\bibitem{CLip}
R. Cluckers, L. Lipshitz, \textit{Fields with analytic structure},
preprint available at http://www.dma.ens.fr/$\sim$cluckers/.









\bibitem{cohen}P. J. Cohen,
\textit{Decision procedures for real and {$p$}-adic fields}, Comm.
Pure Appl. Math. \textbf{22} (1969), 131--151.

\bibitem{Comte}
G. Comte, \textit{Formule de Cauchy-Crofton pour la densit\'e des ensembles sous-analytiques}, C. R.
Acad. Sci. Paris, S\'erie I 328 (1999) 505--508.


\bibitem{Comte2}
G. Comte, \textit{\'Equisingularit\'e r\'eelle : nombres de {L}elong
et images polaires}, Ann. Sci. {\'E}cole Norm. Sup., \textbf{33}
(2000), 757--788.


\bibitem{ComteMerle}
G. Comte, M. Merle,
\textit{\'Equisingularit\'e r\'eelle II :
invariants
locaux et conditions de r\'egularit\'e},
Ann. Sci. {\'E}cole Norm. Sup., \textbf{41}
(2008), 221--269.




\bibitem{D84}
J. Denef, \textit{The rationality of the {P}oincar\'e series
associated to the $p$-adic points on a variety}, Invent. Math.
\textbf{77} (1984), 1--23.

%


\bibitem{Dcell}
J. Denef, \textit{$p$-adic semi-algebraic sets and cell
decomposition}, J. Reine Angew. Math. \textbf{369} (1986), 154--166.



\bibitem{DvdD}
J. Denef, L. van den Dries, \textit{$p$-adic and real subanalytic
sets}, Ann. of Math., \textbf{128}, (1988), 79--138.




\bibitem{vdDriesSkolem}
{L. van den Dries}, \textit{Algebraic theories with definable {S}kolem functions}, The Journal of Symbolic Logic,
\textbf{49} No. 2, (1984), 625--629.


\bibitem{vdDries}
{L. van den Dries}, \textit{Dimension of definable sets, algebraic
boundedness and {H}enselian fields}, Ann. Pure Appl. Logic,
\textbf{45} (1989), 189--209.




\bibitem{DHM} L.~van~den~Dries, D.~Haskell and D.~Macpherson,
\emph{One-dimensional $p$-adic subanalytic sets}, J. London Math.
Soc. (2), \textbf{59}, (1999), 1--20.


\bibitem{sd}
L. van den Dries, P. Scowcroft, \textit{On the structure of
semialgebraic sets over $p$-adic fields}, J. Symbolic Logic,
\textbf{53} (1988), 1138--1164.





\bibitem{Kurdyka}
K. Kurdyka,
\textit{On a subanalytic stratification satisfying a {W}hitney property with exponent $1$.} Real algebraic geometry (Rennes, 1991), 316--322,
Lecture Notes in Math., 1524, Springer, Berlin, 1992.

\bibitem{KR}
K. Kurdyka, G. Raby, \textit{Densit\'e des ensembles
sous-analytiques}, Ann. Inst. Fourier, \textbf{ 39}
(1989), 753--771.



\bibitem{Macint}
A. Macintyre, \textit{On definable subsets of $p$-adic fields},
J. Symbolic Logic \textbf{41} (1976), no. 3, 605--610.

%
%
%


\bibitem{Paru}
A. Parusi\'nski,
\textit{Lipschitz properties of semi-analytic sets},
Ann. Inst. Fourier  \textbf{38} (1988), 189--213.


\bibitem{Pawlucki}
W. Paw{\l}ucki,
 \textit{Lipschitz cell decomposition in $o$-minimal structures. I}, to appear in IIlinois Journal of Mathematics.



\bibitem{Thie}
P. Thie,  \textit{The Lelong number of a point of a complex analytic
set},  Math. Ann.  \textbf{172} (1967) 269--312.

%


\end{thebibliography}

\end{document}